# HARRIS FAMILY OF DISCRETE DISTRIBUTIONS

E, Sandhya; S, Sherly and N, Raju

ABSTRACT. In this paper we discuss the basic properties of a discrete distribution introduced by Harris in 1948 and obtain a characterization of it. The divisibility properties of the distribution are also studied. We derive the moment and maximum likelihood estimators for both the parameters and verify them by simulated observations.

1. INTRODUCTION. The following probability generating function (PGF) was introduced by Harris (Harris, 1948)

$$P(s)=\frac{s}{(m-(m-1)s^k)^{1/k}}, k>0 \text{ integer and } m>1. \quad (1)$$

He discussed this PGF while considering a simple discrete branching process where a particle either splits in to $(k+1)$ identical particles or remains the same during a short time interval $\Delta t$. We refer to the distribution corresponding to the PGF (1) as Harris distribution and denote it by $H_1(m,k,1/k)$. In the notation the suffix 1 suggests that the support of the distribution starts from unity, $m$ determines the probabilities, $k$ implies that the atoms of the distribution are $k$ integers apart and $1/k$ is the shape parameter.

From the PGF one can readily see that this is a generalization of the positive geometric distribution to which it reduces when $k=1$ and its atoms (probability carrying integers) are $k$ integers apart or the probabilities are concentrated on the points $1, 1+k, 1+2k, \ldots$ . It is also true that these probabilities coincide with that of the negative binomial distribution on $\{0, 1, 2, \ldots\}$ with parameters $1/m$ and $1/k$. But no study concerning the properties of this distribution is seen in the literature. In this paper certain basic features of this distribution are brought out.

The main interest that led us to a close examination of this distribution is its role played in schemes with random ($N$) sample sizes (random-sums or $N$-sums and random-extremes or $N$-extremes) in general and in branching processes and time series models in particular where $N$ is a non-negative integer-valued random variable (RV). This PGF had been discussed in the context of branching processes (Harris, 1948), $N$-sums and $N$-extremes (Satheesh *et al.*, 2002, Satheesh and Nair, 2002*a*) and a time series model developed (Satheesh *et al.*, 2005) that has an inherent $N$-sum structure where $N$ is Harris distributed. Also, it is known (Satheesh and Nair, 2004) that a Harris-sum of Harris distributions is again Harris and they used this property to generalize the Marshal-Olkin parametrizing scheme. Harris distribution had been used to demonstrate the notion of random infinite divisibility w.r.t non-negative integer-valued RVs (Sandhya, 1996) and it had also been used to demonstrate that

The research of S. Sherly is supported by the University Grants Commission of India.

an integer-valued infinitely divisible RV $X$ with $P\{X=1\}>0$ can have gaps in its support (Satheesh, 2004*a*). It is also worth mentioning here that in the context of $N$-maximums a reparametrization of this distribution under the name 'extended geometric distribution' had been considered (Voorn, 1987).

As already mentioned $H_1(m,k,1/k)$ denotes the Harris distribution on $\{1, 1+k, 1+2k, \ldots\}$. Other member of the Harris family of distributions considered here is that on $\{0, k, 2k, \ldots\}$ with PGF, $P(s) = (m-(m-1)s^k)^{-1/k}$, $k>0$ integer and $m>1$ and we denote this by $H_0(m,k,1/k)$.

In section.2 we have discussed the genesis and the distributional properties of the Harris family of distributions and a characterization of it. The divisibility properties of the family are studied in section.3. In section.4 we discuss estimation of the parameters and verify them by simulation of the distribution.

2. GENESIS AND DISTRIBUTIONAL PROPERTIES OF HARRIS FAMILY.

(*a*) Let $X$ be a RV degenerate at $k>0$ integer and $Y$ be a negative binomial RV with parameters $1/m$ and $1/k$ ($NB(1/m,1/k)$) having PGF $P(s) = (m-(m-1)s)^{-1/k}$ and let $N = 1+Y$. Then the $N$-sum of independent observations on $X$ is $H_1(m,k,1/k)$ distributed.

(*b*) Consider $\varphi_\theta(s) = (1+\theta s)^{-1/k}$, $k>0$ integer, $\theta>0$, the Laplace transform (LT) of the gamma$(\theta,1/k)$ distribution and let $\varphi(s) = \varphi_1(s)$.

Then $\varphi(\varphi_\theta^{-1}(s)) = \{\frac{1}{\theta}-(\frac{1}{\theta}-1)s^k\}^{-1/k}$, which is the PGF of $H_1(1/\theta,k,1/k)$.

(*c*) Next we obtain the Harris distribution as a gamma mixture of Poisson distribution. Consider the RV $X = kY+1$ where $Y$ has a conditional Poisson distribution with parameter $\lambda$ so that

$$P((Y=r)/\lambda) = \frac{e^{-\lambda}\lambda^r}{r!}, r = 0,1,2, \ldots$$

where $\lambda$ itself is assumed to have a gamma distribution with density

$$g(\lambda) = \frac{\left(\frac{1}{m-1}\right)^{\frac{1}{k}}}{\Gamma(\frac{1}{k})} e^{-\frac{\lambda}{m-1}} \lambda^{\frac{1}{k}-1}, \quad k>0 \text{ integer}, m>1, \lambda>0.$$

The RV $X$ has support $1, 1+k, 1+2k, \ldots$ and the conditional distribution is;

$$P((X=1+rk)/\lambda) = P((Y=r)/\lambda) = \frac{e^{-\lambda}\lambda^r}{r!}, \ r = 0,1,2, \ldots$$

Then the marginal distribution of $X$ is given by

$$P\left(X=1+rk\right) = \int_0^\infty \frac{\left(\frac{1}{m-1}\right)^{\frac{1}{k}}}{\Gamma(\frac{1}{k})} e^{-\frac{\lambda}{m-1}} \lambda^{\frac{1}{k}-1} \frac{e^{-\lambda}\lambda^r}{r!} d\lambda$$

$$= \frac{\left(\frac{1}{m-1}\right)^{\frac{1}{k}}}{r!\ \Gamma(\frac{1}{k})} \int_0^\infty e^{-\left(\frac{1}{m-1}+1\right)\lambda} \lambda^{\frac{1}{k}+r-1} d\lambda$$

$$= \begin{pmatrix} \frac{1}{k}+r-1 \\ r \end{pmatrix} \left(\frac{1}{m}\right)^{\frac{1}{k}} \left(1-\frac{1}{m}\right)^r, \quad r = 0, 1, 2, \ldots$$

using the gamma integral. Thus the marginal distribution of $X$ is $H_1(m,k,1/k)$.

The relationship between the Harris and the negative binomial distribution can also be described as in the following theorem. The proof follows from the corresponding PGFs.

*Theorem*.2.1 A RV $X \sim H_1(m,k,1/k)$ if and only if $Y = \frac{X-1}{k} \sim NB\left(\frac{1}{m},\frac{1}{k}\right)$.

$X \sim H_0(m,k,1/k)$ if and only if $Y = \frac{X}{k} \sim NB\left(\frac{1}{m},\frac{1}{k}\right)$.

*Remark*.2.1 A point worth noticing here is that usually integer-valued RVs are not closed under division. However, here an integer-valued RV is obtained by dividing another integer-valued RV by a constant.

Certain relations between the Harris family, geometric, negative binomial and gamma distributions are given below;

(i) When $k = 1$ in the PGF (1) we get the PGF of the geometric distribution on {1, 2, …}.

$$P(s) = \frac{s}{m-(m-1)s} = \frac{ps}{1-qs},\ p = \frac{1}{m},\ q = 1-p$$

(ii) The PGF $\frac{s}{m-(m-1)s^k} = \frac{ps}{1-qs^k}$ is the PGF of the geometric distribution on {1, 1+$k$, 1+2$k$, …}.

(iii) The PGF $\frac{1}{m-(m-1)s^k} = \frac{p}{1-qs^k}, p = \frac{1}{m}$ is the PGF of the geometric distribution on {0, $k$, 2$k$, …}.

(iv) The PGF (1) can also be written as; $P(s) = \left(\dfrac{s^k}{m-(m-1)s^k}\right)^{\frac{1}{k}} = \left(G(s)\right)^{\frac{1}{k}}$ where

$G(s) = \dfrac{s^k}{m-(m-1)s^k} = \dfrac{ps^k}{1-qs^k}$ is the PGF of the geometric distribution on {$k$, $2k$, $3k$…}. Hence the $k^{th}$ root of the PGF of the geometric distribution on {$k$, $2k$, $3k$…} is the PGF of $H_1(m,k,1/k)$.

(v) Harris PGF (1) is the special case of $\theta = 1$ in the PGF $P^*(s) = \dfrac{s^\theta}{\left(m-(m-1)s^k\right)^{\frac{1}{k}}}$, $\theta>0$ integer, which is the PGF of the linear function $Z = kY + \theta$ of $Y$ where $Y \sim NB\left(\dfrac{1}{m}, \dfrac{1}{k}\right)$.

(vi) The following relation between the PGF $P(s) = (m-(m-1)s)^{-\beta}$ of a $NB\left(\dfrac{1}{m}, \beta\right)$ distribution and the LT $\varphi_\theta(s) = (1+\theta s)^{-\beta}$ of a gamma($\theta,\beta$) is also worth recording here.

$P(s) = \{p/[1-(1-p)s]\}^{-\beta} = \{1+\theta(1-s)\}^{-\beta} = \varphi_\theta(1-s)$, where $p = 1/m$, $q = 1-p$ and $\theta = q/p$.

Now, let $X \sim H_1(m,k,1/k)$. Then from the coefficient of $s^{1+rk}$ in the PGF (1);

$$P(X = 1+rk) = \binom{\frac{1}{k}+r-1}{r}\left(\frac{1}{m}\right)^{\frac{1}{k}}\left(1-\frac{1}{m}\right)^r$$

$$= \binom{-\frac{1}{k}}{r} p^{\frac{1}{k}} (-q)^r, r = 0,1,2, \ldots$$

$k>0$ integer, $m>1$. This is the $(r+1)^{th}$ term in the expansion of $\left(\dfrac{1}{p} - \dfrac{q}{p}\right)^{-\frac{1}{k}}$ *where* $p = \dfrac{1}{m}$, $q = 1-\dfrac{1}{m}$. These probabilities are attached to the points 1, $1+k$, $1+2k$, … which are $k$ integers apart. The recurrence relation for these probabilities is;

$$P\left(X = 1+(r+1)k\right) = \frac{1+rk}{1+r}\,\frac{\left(1-\frac{1}{m}\right)}{k}\,P\left(X = 1+rk\right), r = 0, 1, 2, \ldots$$

The distribution function (DF) is given by;

$$F(1+rk) = \sum_{j=0}^{r}\binom{\frac{1}{k}+j-1}{j}\left(\frac{1}{m}\right)^{\frac{1}{k}}\left(1-\frac{1}{m}\right)^j, r = 0, 1, 2, \ldots$$

This sum can be expressed in terms of an incomplete beta function ratio. Adopting an approach similar to the one in (Johnson *et al.* , 1992) and using theorem.2.1 we have;

$$P\left(X > 1 + (r-1)\,k\right) = P\left(Y > r-1\right)$$

where $X$ is $H_1(m,k,1/k)$ and $Y$ is $NB\left(\frac{1}{m},\frac{1}{k}\right)$. Now proceeding as in (Johnson *et al.* , 1992) we can prove:

$$P\left(Y > r-1\right) = \frac{B_{1-\frac{1}{m}}\left(r,\frac{1}{k}\right)}{B\left(r,\frac{1}{k}\right)} = I_{1-\frac{1}{m}}\left(r,\frac{1}{k}\right) \text{ and hence}$$

$$P\left(X > 1 + rk\right) = I_{1-\frac{1}{m}}\left(r+1,\frac{1}{k}\right).$$

Hence the DF is;

$$F\left(1+rk\right) = 1 - I_{1-\frac{1}{m}}\left(r+1,\frac{1}{k}\right) = I_{\frac{1}{m}}\left(\frac{1}{k}, r+1\right),\; r = 0, 1, 2, \;\cdots$$

where $I_p(a,b) = \dfrac{B_p(a,b)}{B(a,b)}$.

*Theorem*.2.2 Let $X$ and $Y$ be independent and identically distributed (IID) $H_1(m,k,1/k)$ RVs. Then the conditional distribution of $X$ given $X+Y$ is;

$$P\left(X = 1 + rk \,/\, X + Y = 2 + tk\right) = \frac{\binom{\frac{1}{k}+r-1}{r}\binom{\frac{1}{k}+t-r-1}{t-r}}{\binom{\frac{2}{k}+t-1}{t}}$$

$r =$ 0, 1, 2, …, $t$ and $t =$ 0, 1, 2, … . In particular, if $k = 1$, the conditional distribution is uniform.

*Proof.*

$$P\left(X = 1 + rk\right) = P\left(Y = 1 + rk\right) = \binom{\frac{1}{k}+r-1}{r}\left(\frac{1}{m}\right)^{\frac{1}{k}}\left(1-\frac{1}{m}\right)^{r}$$

$$P\left(X = x \,/\, X+Y = z\right) = P\left(X = 1 + rk \,/\, X + Y = 2 + tk\right)$$

$$= \frac{P\left(X = 1 + rk,\; Y = 1 + (t-r)k\right)}{P\left(X + Y = 2 + tk\right)}$$

$$= \frac{P\left(X = 1 + rk\right)P\left(Y = 1 + (t-r)k\right)}{\sum_{r=0}^{t} P\left(X = 1 + rk\right)P\left(Y = 1 + (t-r)k\right)}$$

$$= \frac{\binom{\frac{1}{k} + r - 1}{r}\binom{\frac{1}{k} + t - r - 1}{t - r}}{\sum_{r=0}^{t}\binom{\frac{1}{k} + r - 1}{r}\binom{\frac{1}{k} + t - r - 1}{t - r}}$$

$$= \frac{\binom{\frac{1}{k} + r - 1}{r}\binom{\frac{1}{k} + t - r - 1}{t - r}}{\binom{\frac{2}{k} + t - 1}{t}}$$

If $k = 1$, that is, if $X$ and $Y$ are IID geometric RVs, then

$$P\left(X = 1 + r \,/\, X + Y = 2 + t\right) = \frac{1}{1+t},\ r = 0, 1, 2, \ldots t;\ t = 0, 1, 2, \cdots$$

and hence the conditional distribution is uniform.

Let $\mu_{(r)}$, $\mu'_r$ and $\mu_r$ denote respectively the r<sup>th</sup> factorial moment, r<sup>th</sup> moment about zero and r<sup>th</sup> central moment of a distribution. We now evaluate the first four of them for the $H_1(m,k,1/k)$ distribution. From this one can identify the mean and variance of the $H_1(m,k,1/k)$ distribution.

$$\mu_{(1)} = m$$
$$\mu_{(2)} = m\,(m-1)\,(k+1)$$
$$\mu_{(3)} = m\,(m-1)\,(k+1)\,\left(k - m + 2\,(k+1)\,(m-1)\right)$$
$$\mu_{(4)} = m\,(m-1)\,(k+1)\,\left(6m^2k^2 + 5m^2k + m^2 - 6mk^2 - 13mk - 5m + k^2 + 5k + 6\right)$$

$$\mu'_1 = m$$

$$\mu'_2 = m^2 + m\,(m-1)\,k$$

$$\mu'_3 = m^3 + m\,(m-1)\,k\,\left(2mk + 3m - k\right)$$

$$\mu'_4 = m^4 + m\,(m-1)\,k\,\left(6m^2k^2 + 11m^2k + 6m^2 - 6mk^2 - 7mk + k^2\right)$$

Similarly,

$$\mu_1 = 0$$

$$\mu_2 = m\,(m-1)\,k$$

$$\mu_3 = m\,(m-1)\,(2m-1)\,k^2$$

$$\mu_4 = m\,(m-1)\,\big(k(6m^2-6m+1) + 3\,m\,(m-1)\big)\,k^2.$$

Harris distribution is positively skewed since;

$$\beta_1 = \frac{(2m-1)^2\;k}{m\;(m-1)}$$

$$\gamma_1 = (2m-1)\sqrt{\frac{k}{m\,(m-1)}}$$

and $\mu_3 = m\,(m-1)\,(2m-1)\;k^2 > 0$.

Harris distribution is leptokurtic as $\beta_2 > 3$;

$$\beta_2 = 3 + 6k + \frac{k}{m\,(m-1)}$$

$$\gamma_2 = 6k + \frac{k}{m\,(m-1)}$$

It's coefficient of variation is; C.V $= \sqrt{\left(1 - \frac{1}{m}\right)k}$ .

From the PGF of $H_1(m,k,1/k)$ it follows that the moment generating function is;

$$\frac{e^t}{(m-(m-1)\,e^{tk})^{\frac{1}{k}}}$$

and so the cumulant generating function is;

$$t - \frac{1}{k}\log(m) - \frac{1}{k}\log\left(1-\left(1-\frac{1}{m}\right)e^{tk}\right).$$

Hence the first four cumulants are;

$$k_1 = m$$

$$k_2 = m\,(m-1)\,k$$

$$k_3 = m\,(m-1)\,(2m-1)\,k^2$$

$$k_4 = m\,(m-1)\,(6m^2-6m+1)\,k^3$$

Now, the recurrence relations for moments are

$$\mu'_r = \sum_{j=0}^{\infty} (1+jk)^r \, P(X=1+jk)$$ gives the relationship

$$\mu'_{r+1} = m\left( (m-1)\,k\,\frac{d\mu'_r}{d\,m} + \mu'_r \right) \quad , r = 0,\ 1,\ 2,\ \ldots$$

$$\mu_r = \sum_{j=0}^{\infty} (1+jk-m)^r \, P(X=1+jk)$$ gives

$$\mu_{r+1} = m\,(m-1)\,k\left( \frac{d\mu_r}{dm} + r\,\mu_{r-1} \right) \quad , r = 1,\ 2,\ \ldots$$

The recurrence relation for cumulants ;

Since the r$^{th}$ cumulant $k_r = \left[ \frac{d^r}{dt^r} K_X(t) \right]_{t=0}$ and $\frac{dk_r}{dm} = \left[ \frac{d^r}{dt^r} \frac{d}{dm} K_X(t) \right]_{t=0}$

the cumulants satisfy the relationship $k_{r+1} = m\,(m-1)\,k\,\frac{dk_r}{dm}$ , $r = 1,\ 2,\ \ldots$

From the relation

$$\frac{P(X=1+(r+1)k)}{P(X=1+rk)} = \frac{1+rk}{k+rk}\left(1-\frac{1}{m}\right), \quad r = 0,\ 1,\ 2,\ \ldots\ ,$$

it can be seen that;

$$P(X=1+rk) > P(X=1+(r+1)k), \quad r = 0,\ 1,\ 2,\ \ldots$$

and there is a single mode at $X = 1$.

Now we present a characterization of the Harris distribution motivated by a characterization of the negative binomial distribution (Sathe and Ravi, 1997).

*Theorem*.2.3 Let $X$ be a RV with $P\{X = nk+1\} = a_{nk+1}$, n = 0,1,2, … and let $E(X) = \mu$ be finite.

Then $X \sim H_1(m,k,1/k)$ if and only if

$$\frac{dP(X > nk+1)}{d\mu} = \frac{nk+1}{\mu\,k}\, a_{nk+1}, \quad n = 0,1,2,\ \ldots,\ k>0 \text{ is an integer} \qquad (2)$$

*Proof.* Let $X \sim H_1(m,k,1/k)$ and let $E(X) = \mu$ be finite.

$$P(X > nk+1) = 1 - \sum_{r=0}^{n} \binom{\frac{1}{k} + r - 1}{r} \left(\frac{1}{m}\right)^{\frac{1}{k}} \left(1 - \frac{1}{m}\right)^{r}$$

$$\frac{d}{d\mu} P(X > nk+1) = -\frac{d}{d\mu} \sum_{r=0}^{n} \binom{\frac{1}{k} + r - 1}{r} \left(\frac{1}{\mu}\right)^{\frac{1}{k}} \left(1 - \frac{1}{\mu}\right)^{r}$$

Since for $H_1(m,k,1/k)$, $\mu = E(X) = m$. Differentiation and simplification gives condition (2). Assuming the condition (2) and using $P(X > nk+1) = 1 - P(X \le nk+1)$ we have

$$-\sum_{r=0}^{n} \frac{d}{d\mu} a_{rk+1} = \frac{nk+1}{\mu k} a_{nk+1}, \; n = 0, 1, 2, \ldots$$

Consequently , we get

$$-\frac{da_1}{d\mu} = \frac{a_1}{\mu k} \quad \text{and}$$

$$-\frac{d\, a_{nk+1}}{d\mu} = \frac{nk+1}{\mu k} a_{nk+1} - \frac{(n-1)k+1}{\mu k} a_{(n-1)k+1}, \quad n \ge 1 \tag{3}$$

Defining $P(\mu, s) = \sum_{n=0}^{\infty} a_{nk+1} s^{nk+1}$, $0 < s \le 1$, we get

$$\frac{\partial P}{\partial s} = \sum_{n=0}^{\infty} (nk+1) a_{nk+1} s^{nk} \quad \text{and} \quad \frac{\partial P}{\partial \mu} = s\frac{\partial a_1}{\partial \mu} + \sum_{n=1}^{\infty} \frac{\partial a_{nk+1}}{\partial \mu} s^{nk+1} .$$

Thus from (3) we get $-\mu k \dfrac{\partial P}{\partial \mu} = s(1 - s^k)\dfrac{\partial P}{\partial s}$

Consequently, we get the system of differential equations

$$\frac{d\mu}{\mu k} = \frac{ds}{s(1-s^k)} = \frac{dp}{0}$$

Solving these equations with the initial conditions $\mu = 1$ and $P(1, s) = s$ we get

$$P(\mu, s) = \left(\frac{s^k}{\mu - \mu s^k + s^k}\right)^{\frac{1}{k}}$$

which is the PGF of $H_1(m,k,1/k)$ with probability mass function $p_{nk+1}$. This completes the proof.

The following figures give us an idea about the shape of the Harris distribution. They correspond to the probability distribution given below.

| Harris (50,5,1/5) | |
|---|---|
| x | P(x) |
| 1 | 0.457305 |
| 6 | 0.089632 |
| 11 | 0.052704 |
| 16 | 0.037876 |
| 21 | 0.029695 |
| 26 | 0.024445 |
| 31 | 0.020762 |
| 36 | 0.018021 |
| 41 | 0.015895 |
| 46 | 0.014192 |
| 51 | 0.012796 |
| 56 | 0.011628 |
| 61 | 0.010636 |
| 66 | 0.009781 |
| 71 | 0.009038 |
| 76 | 0.008385 |
| 81 | 0.007806 |
| 86 | 0.00729 |
| 91 | 0.006827 |
| 96 | 0.006409 |
| 101 | 0.006029 |

| Harris (4,5,1/5) | |
|---|---|
| x | P(x) |
| 1 | 0.757858 |
| 6 | 0.113679 |
| 11 | 0.051155 |
| 16 | 0.028136 |
| 21 | 0.016881 |
| 26 | 0.010635 |
| 31 | 0.006913 |
| 36 | 0.004592 |
| 41 | 0.0031 |
| 46 | 0.002118 |
| 51 | 0.001462 |
| 56 | 0.001016 |
| 61 | 0.000711 |
| 66 | 0.000501 |
| 71 | 0.000354 |
| 76 | 0.000251 |
| 81 | 0.000179 |
| 86 | 0.000128 |
| 91 | 9.18E-05 |
| 96 | 6.59E-05 |

| Harris (2,2,1/2) | |
|---|---|
| x | p(x) |
| 1 | 0.707107 |
| 3 | 0.176777 |
| 5 | 0.066291 |
| 7 | 0.027621 |
| 9 | 0.012084 |
| 11 | 0.005438 |
| 13 | 0.002492 |
| 15 | 0.001157 |
| 17 | 0.000542 |
| 19 | 0.000256 |
| 21 | 0.000122 |
| 23 | 5.81E-05 |
| 25 | 2.78E-05 |
| 27 | 1.34E-05 |

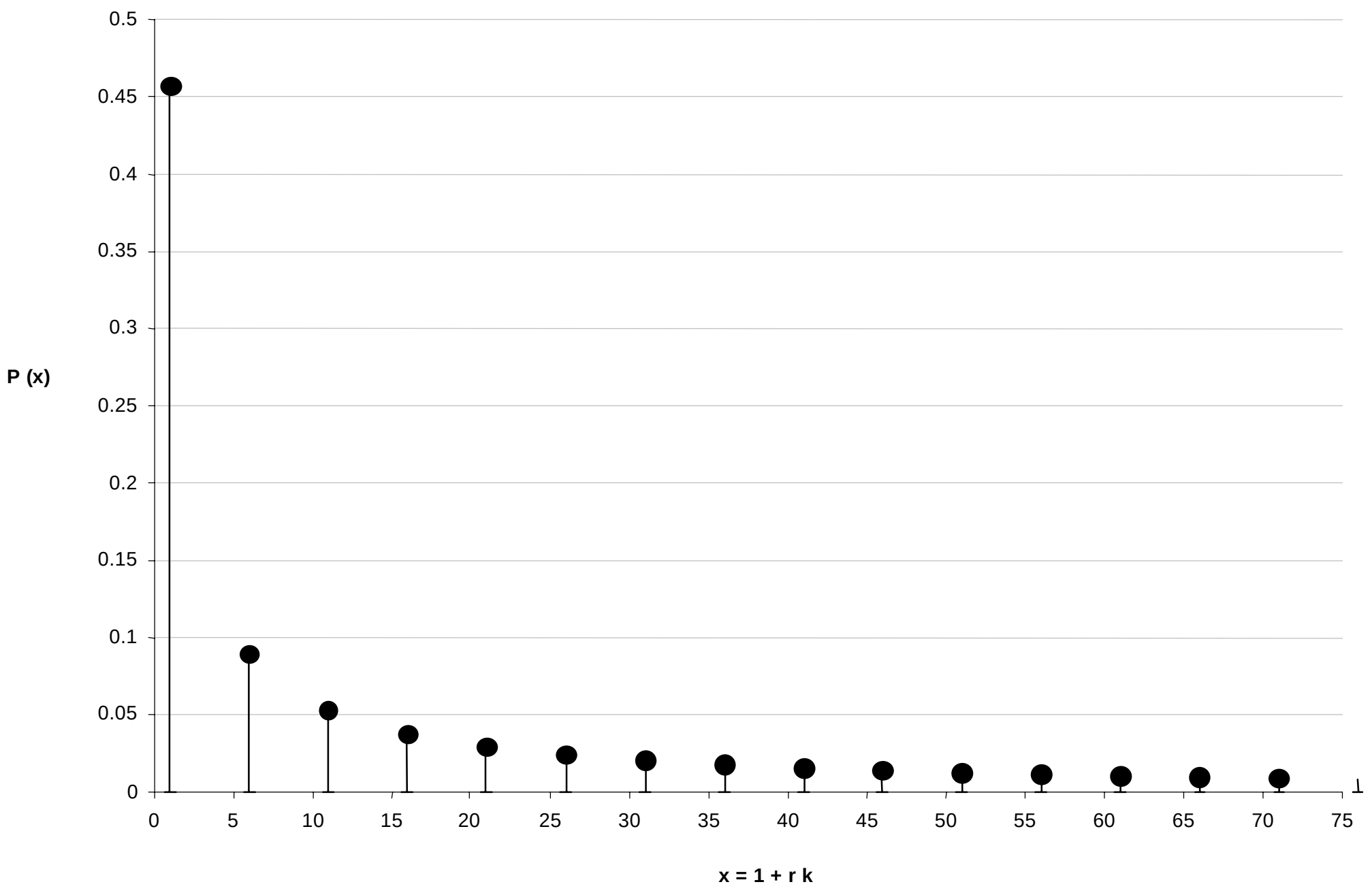

*Figure 1* - Probability Plot of Harris(50,5,1/5)

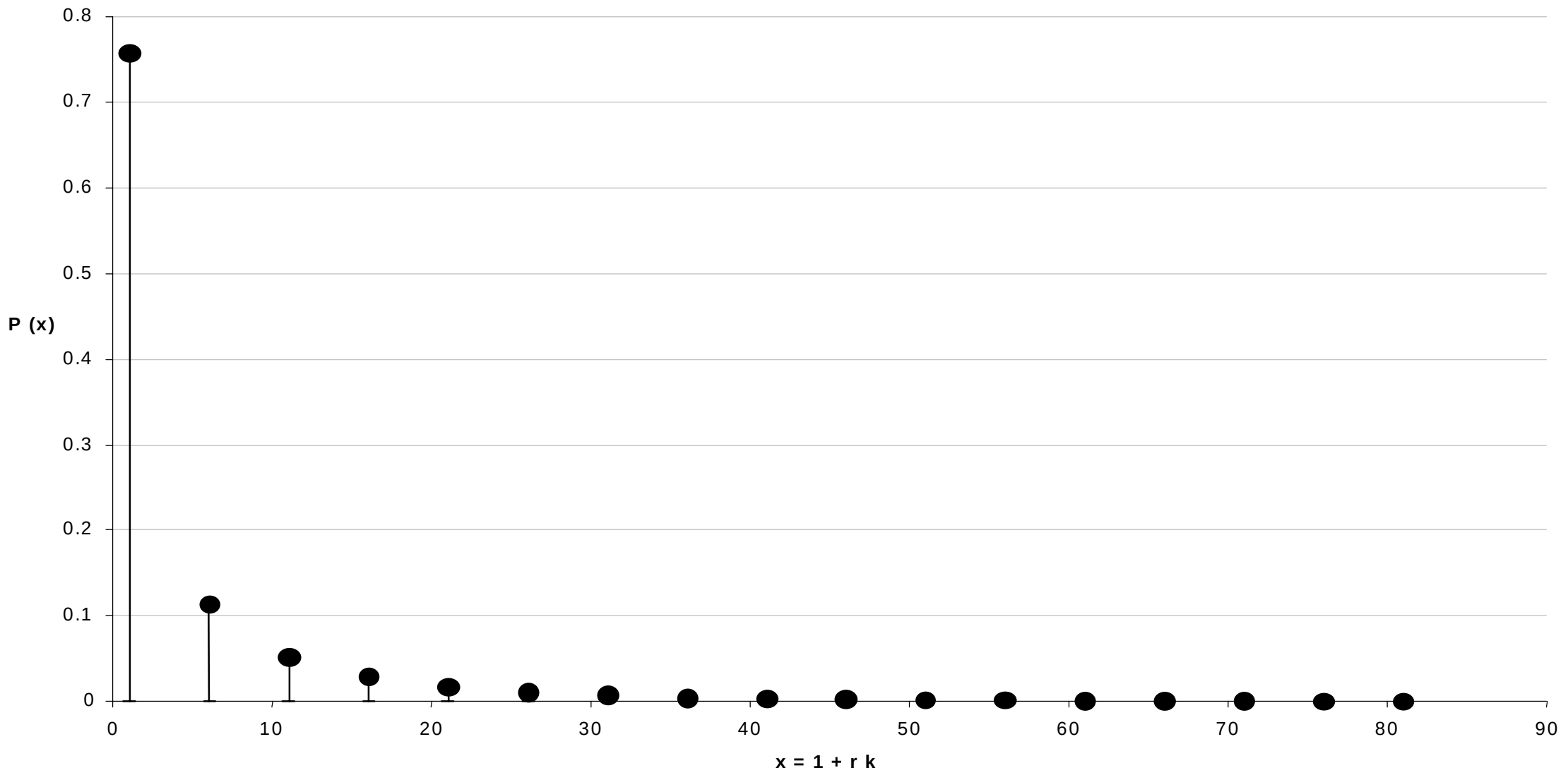

*Figure 2* – Probability Plot of Harris(4, 5, 1/5).

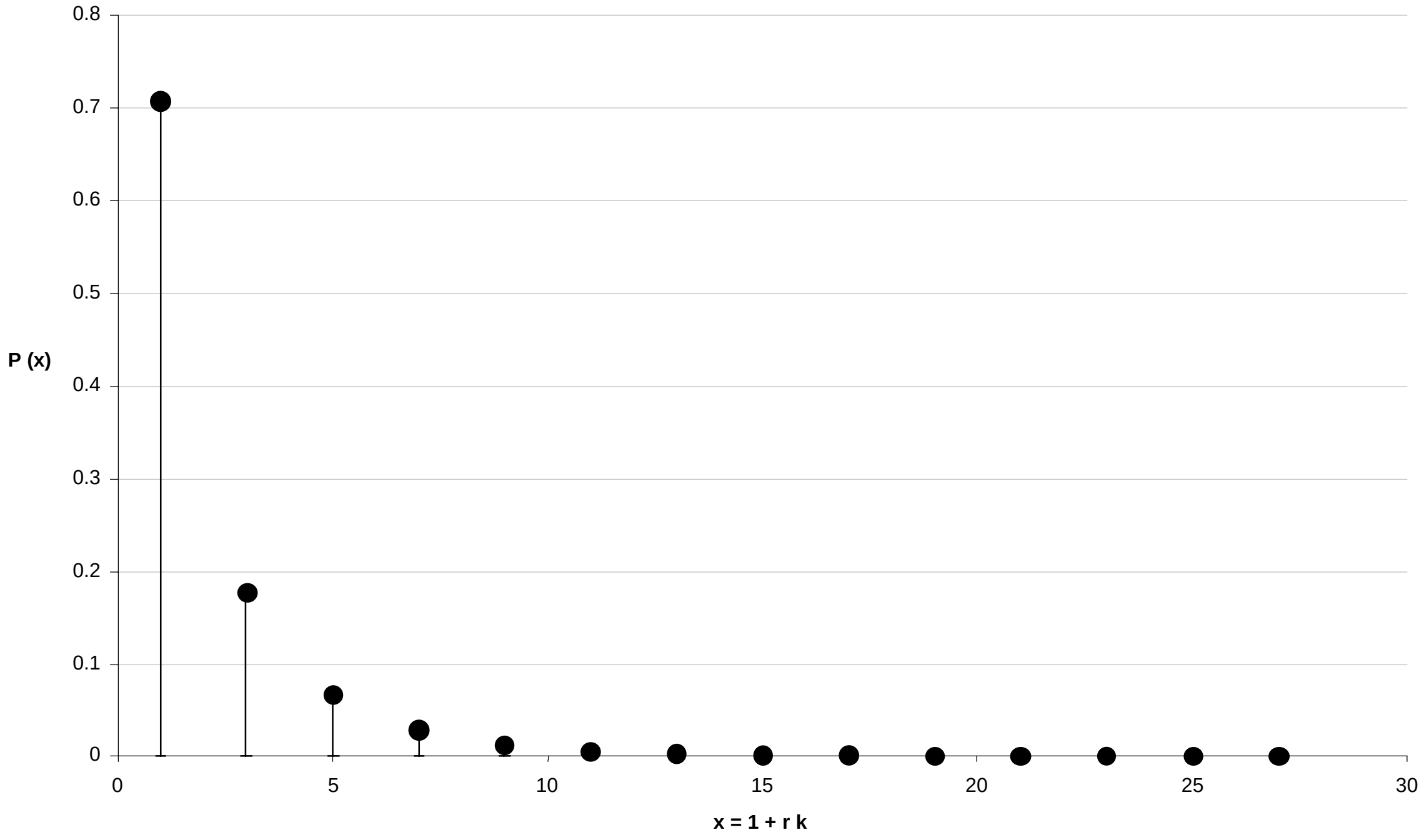

*Figure 3* – Probabilty Plot of Harris(2, 2, ½).

## 3. DIVISIBILITY PROPERTIES OF HARRIS FAMILY.

The PGF of the $H_0(m,k,1/k)$ distribution is $P(s) = \{m-(m-1)s^k\}^{-1/k}$ which corresponds to a random variable $U = kV$, where $V$ is $NB(1/m,1/k)$ with PGF $\{m-(m-1)s\}^{-1/k}$. We know that a negative binomial distribution is infinitely divisible (ID). Since infinite divisibility is not affected by change of origin and scale, it follows that $U$ is ID. Hence;

*Theorem*.3.1 $H_0(m,k,1/k)$ and $H_1(m,k,1/k)$ RVs are ID .

We know that gamma($c,\beta$) distributions are self-decomposable (Sreehari, 1979). We now give an alternate proof that gamma($c,1/k$) distributions are self-decomposable using its Harris–sum stability. As a consequence of this we get that $H_0(m,k,1/k)$ distributions are self-decomposable. (Satheesh *et al.* ,2002) showed that gamma($c,1/k$) law is $H_1(a,k,1/k)$-sum stable in the sense:

$$\frac{1/(1+cs)^{\frac{1}{k}}}{[a-(a-1)/(1+cs)]^{\frac{1}{k}}} = \frac{(1+cs)^{-\frac{1}{k}}}{[(1+acs)/(1+cs)]^{\frac{1}{k}}} = \left(\frac{1}{1+acs}\right)^{\frac{1}{k}}.$$

From the above when $a = 1/c$ we have:

$$\left(\frac{1}{1+s}\right)^{\frac{1}{k}} = \left(\frac{1}{1+cs}\right)^{\frac{1}{k}} \frac{1}{\left[a-(a-1)\dfrac{1}{1+cs}\right]^{\frac{1}{k}}}.$$

Here the second factor on the right hand side is a LT being that of a $H_0(a,k,1/k)$-sum of gamma($c,1/k$) variables. The relation is also true for each $0<c<1$ since we may choose the $H_0(a,k,1/k)$ with $a = 1/c$ accordingly. Hence gamma($c,1/k$) RV is self-decomposable. Consequently the corresponding negative binomial RVs are also self-decomposable by corollary.2.1 in (Satheesh and Nair, 2002*a*). Since self-decomposability is not affected by change of scale we have proved;

*Theorem*.3.2 $H_0(a,k,1/k)$ RVs are self-decomposable.

*Remark*.3.1 $H_1(a,k,1/k)$ RVs are not self-decomposable as it has no probability at zero, a necessary condition for a discrete distribution to be self-decomposable, (Satheesh, 2004*a*).

(Kozubowski and Panorska, 1996) had developed ν-stable laws and according to them the ν-stable law with characteristic function (CF ) $\{1+|t|^{\alpha}\}^{-\beta}$ is obtained as the weak limit of negative binomial sums. The above CF is that of the generalized Linnik distribution and gamma($c,\beta$) distribution is a particular case of it. But negative binomial RV does not induce $N$-sum stability in the generalized Linnik distributions. It is the $H_1(a,k,1/k)$ RV that induces $N$-sum stability in generalized Linnik distributions and that too only when $\beta = 1/k$, see (Satheesh *et al.* , 2002). Next we show that the above CF can be obtained as the weak limit of Harris-sums also. Let $N_a \sim H_1(a,k,1/k)$. The next lemma follows from (Satheesh, 2004*b*).

*Lemma*.3.1 $N_a \xrightarrow{d} kU$ as $a \to \infty$ where $U$ is gamma$(c, 1/k)$, $c = 1/a$.

*Lemma*.3.2 $N_a \xrightarrow{p} \infty$ as $a \to \infty$.

*Proof.* Let $p_{a,k} = P\{N_a = k\}$, $k = 0, 1, 2, \ldots$ be the probability distribution corresponding to the RVs $\{N_a, a>1\}$. Then $N_a \xrightarrow{p} \infty$ as $a \to \infty$ is equivalent to $p_{a,k} \to 0$ as $a \to \infty$ for every $k = 0, 1, 2, \ldots$ . In terms of PGFs this is equivalent to $\underset{a\to\infty}{Lt}\, P_a(s) = 0$ for all $s \in (0,1)$. Since gamma distribution with LT $\varphi$ is absolutely continuous $\varphi(\infty) = 0$ and since the PGF of $N_a$ is derived from the LT $\varphi$ (see (vi) on page.4), $\underset{a\to\infty}{Lt}\, P_a(s) = 0$. That completes the proof of the lemma.

This lemma enables one to invoke the (Szasz, 1972) transfer theorem where the conditions to be satisfied by the indexing RV $N_a$ are $N_a \xrightarrow{d} \nu$ and $N_a \xrightarrow{p} \infty$ as $a \to \infty$, where $\nu$ is the RV w.r.t which $\nu$-stable laws are described. Here it is that of the gamma$(c,\beta)$ RV. The following result now follows by invoking (Szasz, 1972) theorem.

*Theorem*.3.3 Let $S_n = X_1+X_2+ \ldots + X_n$ where $X_1, X_2, \ldots$are IID RVs and $S_n \xrightarrow{d} U$ so that $U$ is ID with CF $\exp\{-\psi(t)$. Then $S_{N_a} \xrightarrow{d} V$ and the CF of $V$ is $\{1+\psi(t)\}^{-1/k}$ and conversely.

Notice that we may have $N_a \sim H_0(a,k,1/k)$ or $H_1(a,k,1/k)$ and get the result since the conclusion in lemma.3.1 & 3.2 are true for both, see (Satheesh, 2004*b*). Notice also that instead of $N_a \sim H_1(a,k,1/k)$ we may consider the general case of $N_a \sim H_1(a,k,\beta)$ and see that the above arguments hold good. These results extend $\nu$-stable distributions of (Kozubowski and Panorska, 1996) and $\varphi$-stable distributions of (Satheesh, 2004*b*).

## 4. ESTIMATION AND SIMULATION

In the sequel we assume that both the parameters are unknown. We develop maximum likelihood estimators and moment estimators for both parameters $m$ and $k$ of the Harris family of distributions and then evaluate their performances by simulation of the distribution.

*4.1 Method of Maximum Likelihood.* Following the development in (Simon, 1961) we have;

$$f(x) = \binom{\frac{1}{k}+r-1}{r}\left(\frac{1}{m}\right)^{\frac{1}{k}}\left(1-\frac{1}{m}\right)^{r}, \quad \text{where } x = 1+rk, \quad r = 0, 1, 2, \ldots$$

Taking $K = 1/k$ and $p = 1/m$ we have the likelihood function

$$L = \prod_{i=1}^{n}\binom{K+r_i-1}{r_i} p^{nK}(1-p)^{n\bar{r}}, \quad \bar{x} = 1+\bar{r}k$$

$$\log L = \sum_{i=1}^{n} \log \binom{K + r_i - 1}{r_i} + nK \log p + n\bar{r} \log(1-p)$$

$$\frac{\partial \log L}{\partial p} = \frac{nK}{p} - \frac{n\bar{r}}{1-p} = 0 \qquad \text{gives}$$

$$\hat{p} = \frac{K}{K+\bar{r}} = \frac{1}{1+\bar{r}k} = \frac{1}{\bar{\bar{x}}} \qquad \text{and} \quad \text{hence} \quad \hat{m} = \bar{x}$$

The maximum likelihood equation corresponding to $K$ is

$$\sum_{i=1}^{n}\left( \frac{1}{K} + \frac{1}{K+1} + \cdots + \frac{1}{K+r_i-2} + \frac{1}{K+r_i-1} \right) + n \log p + n(\bar{x}-1) \log (1-p) = 0$$

Substituting the value of $p$, we have

$$\sum_{i=1}^{n}\left( \frac{1}{K} + \frac{1}{K+1} + \cdots + \frac{1}{Kx_i-2} + \frac{1}{Kx_i-1} \right) + n \log\left( \frac{1}{\bar{x}} \right) + n\,(\bar{x}-1) \log\left( 1- \frac{1}{\bar{x}} \right) = 0$$

This equation has only one unknown $K$ and , therefore, can be solved by trial and error method.

Then $\hat{k} = 1/\hat{K}$.

*4.2 Method of moments.* Using the method of moments, that is, by equating the sample mean $\bar{x}$ and the sample variance $s^2$ to the corresponding population values we have the following estimates.

$$\hat{m} = \bar{x} \quad \text{and}$$

$$\hat{k} = \frac{s^2}{\bar{x}(\bar{x}-1)}.$$

The results based on simulation are presented in the tables appended.

ACKNOWLEDGEMENT.

The authors thank Dr. T. M. Jacob, Department of Statistics, Nirmala College, Moovattupuzha, for the help in doing the simulation study.

*Department of Statistics* EDAKKUNNY SANDHYA

*Prajyoti Niketan College, Pudukkad, Trichur – 680 301, India.*

*Department of Statistics* SEBASTIAN SHERLY

*Vimala College, Trichur – 680 009, India.*

*Department of Statistics* NARAYANAN RAJU

*University of Calicut, Calicut University PO – 673 635, India.*

## REFERENCES.

TABLE 1

*Maximum likelihood estimates for $k = 2$ and $m = 10$ using simulated samples of different sizes $n$.*

| $n$ | 100 | 200 | 300 | 400 | 500 |
|---|---|---|---|---|---|
| $\hat{k}$ | 2.44 | 2.53 | 2.49 | 2.49 | 2.48 |
| $\hat{m}$ | 8.56 | 10.78 | 10.27 | 9.77 | 9.83 |

TABLE 2

*Maximum likelihood estimates for $k = 2$ and $m = 2$ using simulated samples of different sizes $n$.*

| $n$ | 100 | 200 | 300 | 400 | 500 |
|---|---|---|---|---|---|
| $\hat{k}$ | 2 | 1.91 | 1.92 | 1.98 | 2 |
| $\hat{m}$ | 2.06 | 1.94 | 1.95 | 1.98 | 2 |

TABLE 3

*Maximum likelihood estimates for $k = 4$ and $m = 2$ using simulated samples of different sizes $n$.*

| $n$ | 100 | 200 | 300 | 400 | 500 |
|---|---|---|---|---|---|
| $\hat{k}$ | 4.95 | 5 | 4.81 | 4.83 | 4.95 |
| $\hat{m}$ | 1.72 | 1.78 | 1.91 | 1.89 | 1.85 |

TABLE 4

*Moment estimates of m and k using simulated sample of size 200 and number of repetitions 50.*

| | $k$ | 2 | | 4 | | 10 | | 20 | | 30 | | 50 | |
|---|---|---|---|---|---|---|---|---|---|---|---|---|---|
| $m$ | | $\hat{m}$ | $\hat{k}$ | $\hat{m}$ | $\hat{k}$ | $\hat{m}$ | $\hat{k}$ | $\hat{m}$ | $\hat{k}$ | $\hat{m}$ | $\hat{k}$ | $\hat{m}$ | $\hat{k}$ |
| 1.25 | Estimate | 1.24560 | 2.00013 | 1.27040 | 3.87898 | 1.24300 | 9.51323 | 1.26800 | | 1.23100 | | 1.24500 | |
| | SE | 0.00813 | 0.04954 | 0.01018 | 0.09065 | 0.01732 | 0.30052 | 0.02610 | | 0.03487 | | 0.04255 | |
| 1.5 | Estimate | 1.51660 | 1.96823 | 1.51720 | 3.89563 | 1.51900 | 9.71559 | 1.51000 | 18.88118 | 1.47100 | | 1.49500 | |
| | SE | 0.01307 | 0.04971 | 0.01806 | 0.15076 | 0.02587 | 0.41215 | 0.03695 | 0.84854 | 0.04371 | | 0.05691 | |
| 2 | Estimate | 2.01360 | 1.96947 | 2.03960 | 3.74475 | 1.95000 | 9.13116 | 2.03000 | 17.41241 | 1.87600 | | 1.98000 | |
| | SE | 0.01856 | 0.04596 | 0.02578 | 0.07293 | 0.03774 | 0.29010 | 0.05836 | 0.68539 | 0.08314 | | 0.10567 | |
| 10 | Estimate | 9.78780 | 1.95334 | 9.98520 | 4.01341 | 10.45500 | 9.25310 | 10.25000 | 18.11888 | 9.69500 | 42.60791 | 10.11100 | 28.04847 |
| | SE | 0.13499 | 0.04532 | 0.20464 | 0.09877 | 0.33099 | 0.34030 | 0.48299 | 1.29898 | 0.62402 | 2.79929 | 0.59225 | 1.45120 |
| 50 | Estimate | 50.74180 | 2.02585 | 50.06560 | 4.02037 | 49.14100 | 9.24440 | 50.80200 | 18.94129 | 50.09800 | 27.27051 | 46.09500 | 41.99895 |
| | SE | 0.67331 | 0.06426 | 0.98037 | 0.09138 | 1.24350 | 0.35530 | 2.08304 | 0.86402 | 2.82813 | 1.42751 | 2.80317 | 2.52180 |

TABLE 5

*Moment estimates of m and k using simulated sample of size 100 and number of repetitions 50.*

| | $k$ | 2 | | 4 | | 10 | | 20 | | 30 | | 50 | |
|---|---|---|---|---|---|---|---|---|---|---|---|---|---|
| $m$ | | $\hat{m}$ | $\hat{k}$ | $\hat{m}$ | $\hat{k}$ | $\hat{m}$ | $\hat{k}$ | $\hat{m}$ | $\hat{k}$ | $\hat{m}$ | $\hat{k}$ | $\hat{m}$ | $\hat{k}$ |
| 1.25 | Estimate | 1.26400 | 1.88702 | 1.24320 | 3.73081 | 1.27000 | | 1.24400 | | 1.25800 | | 1.22000 | |
| | SE | 0.01179 | 0.05516 | 0.01498 | 0.10512 | 0.02527 | | 0.03846 | | 0.03737 | | 0.04772 | |
| 1.5 | Estimate | 1.52520 | 1.94526 | 1.50400 | 3.89248 | 1.47600 | 8.55077 | 1.49600 | | 1.46800 | | 1.44000 | |
| | SE | 0.01959 | 0.06138 | 0.02470 | 0.14699 | 0.03839 | 0.39581 | 0.05488 | | 0.05565 | | 0.07778 | |
| 2 | Estimate | 1.96600 | 1.80501 | 2.01600 | 3.79882 | 1.99800 | 8.86396 | 1.82000 | | 1.98400 | | 1.97000 | |
| | SE | 0.02937 | 0.04348 | 0.03977 | 0.13584 | 0.07174 | 0.35652 | 0.09129 | | 0.09073 | | 0.13356 | |
| 10 | Estimate | 9.66800 | 1.87155 | 9.65760 | 4.16408 | 9.66600 | 8.51517 | 10.23200 | 18.46166 | 9.96400 | 19.94747 | 10.57000 | 35.30626 |
| | SE | 0.15140 | 0.06181 | 0.27180 | 0.22783 | 0.40130 | 0.34897 | 0.75848 | 1.49267 | 0.68864 | 1.04824 | 1.03263 | 2.22629 |
| 50 | Estimate | 49.24520 | 1.98421 | 49.06000 | 3.64410 | 45.46400 | 8.73130 | 50.96800 | 16.60633 | 45.85000 | 22.21521 | 50.90000 | 32.10199 |
| | SE | 0.89951 | 0.07132 | 1.24629 | 0.13499 | 1.90484 | 0.39150 | 3.59164 | 0.90848 | 3.32435 | 1.22744 | 5.15364 | 1.99296 |

TABLE 6

*Moment estimates of m and k using simulated sample of size 50 and number of repetitions 100.*

| | $k$ | 2 | | 4 | | 10 | | 20 | | 30 | | 50 | |
|---|---|---|---|---|---|---|---|---|---|---|---|---|---|
| $m$ | | $\hat{m}$ | $\hat{k}$ | $\hat{m}$ | $\hat{k}$ | $\hat{m}$ | $\hat{k}$ | $\hat{m}$ | $\hat{k}$ | $\hat{m}$ | $\hat{k}$ | $\hat{m}$ | $\hat{k}$ |
| 1.25 | Estimate | 1.25040 | 1.91150 | 1.24240 | | 1.23400 | | 1.31200 | | 1.21600 | | 1.27000 | |
| | SE | 0.01035 | 0.04978 | 0.01614 | | 0.02345 | | 0.04083 | | 0.04136 | | 0.05478 | |
| 1.5 | Estimate | 1.52360 | 1.86699 | 1.48640 | 3.85572 | 1.57200 | | 1.46400 | | 1.45600 | | 1.53000 | |
| | SE | 0.01915 | 0.04944 | 0.02459 | 0.11409 | 0.04000 | | 0.05325 | | 0.05665 | | 0.08097 | |
| 2 | Estimate | 2.04920 | 1.86178 | 2.07600 | 3.71790 | 2.06200 | | 2.08400 | | 1.91200 | | 2.07000 | |
| | SE | 0.02837 | 0.05047 | 0.04406 | 0.12435 | 0.06374 | | 0.09102 | | 0.10108 | | 0.15126 | |
| 10 | Estimate | 10.35040 | 1.86544 | 10.35680 | 3.42257 | 10.45400 | 8.62218 | 9.81200 | 14.68714 | 10.17400 | 17.36004 | 8.97000 | 35.30626 |
| | SE | 0.24092 | 0.05620 | 0.25330 | 0.11688 | 0.39496 | 0.44564 | 0.52297 | 0.63116 | 0.75536 | 0.75087 | 1.19405 | 2.22629 |
| 50 | Estimate | 50.66960 | 1.95460 | 51.55760 | 3.83184 | 49.37200 | 7.94068 | 48.45600 | 13.73580 | 53.96200 | 17.85134 | 53.19000 | 32.10199 |
| | SE | 0.92340 | 0.06268 | 1.35669 | 0.17002 | 2.33855 | 0.28221 | 3.21871 | 0.57696 | 3.80399 | 0.72906 | 4.97949 | 1.99296 |

TABLE 7

*Moment estimates of m and k using simulated sample of size 500 and number of repetitions 50.*

| | $k$ | 2 | | 4 | | 10 | | 20 | | 30 | | 50 | |
|---|---|---|---|---|---|---|---|---|---|---|---|---|---|
| $m$ | | Estimate | SE | Estimate | SE | Estimate | SE | Estimate | SE | Estimate | SE | Estimate | SE |
| 1.25 | $\hat{k}$ | 1.99077 | 0.03264 | 4.00879 | 0.07218 | 9.58442 | 0.21471 | 19.17520 | 0.50422 | | | | |
| | $\hat{m}$ | 1.25008 | 0.00479 | 1.24864 | 0.00655 | 1.25800 | 0.01108 | 1.22720 | 0.01483 | | | | |
| 1.5 | $\hat{k}$ | 2.00606 | 0.02766 | 3.92385 | 0.08145 | 9.71578 | 0.22164 | 18.77220 | 0.54252 | 29.1243 | 1.7017 | 48.5713 | 2.8846 |
| | $\hat{m}$ | 1.49208 | 0.00689 | 1.47872 | 0.01080 | 1.50320 | 0.01501 | 1.50320 | 0.02351 | 1.4296 | 0.0250 | 1.5700 | 0.0372 |
| 2 | $\hat{k}$ | 1.98755 | 0.03210 | 3.88569 | 0.07429 | 9.58715 | 0.26509 | 20.55380 | 0.88114 | 25.6596 | 0.9537 | 44.1528 | 1.9316 |
| | $\hat{m}$ | 1.97968 | 0.01220 | 2.00176 | 0.01896 | 1.99800 | 0.02317 | 1.99999 | 0.03869 | 1.9792 | 0.0531 | 1.9940 | 0.0630 |
| 10 | $\hat{k}$ | 1.99643 | 0.03172 | 4.09977 | 0.08131 | 9.44012 | 0.28046 | 19.45160 | 0.60513 | 29.2936 | 1.2349 | 44.8031 | 2.0597 |
| | $\hat{m}$ | 10.02180 | 0.08260 | 10.09070 | 0.10855 | 9.87000 | 0.17870 | 10.37874 | 0.24026 | 9.6772 | 0.3157 | 9.4640 | 0.3755 |
| 50 | $\hat{k}$ | 1.99377 | 0.02259 | 4.06961 | 0.07325 | 9.47269 | 0.29482 | 20.11050 | 0.70908 | 29.2949 | 1.1526 | 44.6115 | 2.0703 |
| | $\hat{m}$ | 50.20820 | 0.51673 | 49.36580 | 0.56156 | 50.9 7960 | 1.14684 | 50.06624 | 1.46995 | 53.0200 | 1.6164 | 50.4960 | 2.1729 |